\documentclass[reqno]{amsproc}

\usepackage{amsmath,amssymb,amsthm,amsfonts,latexsym}
\usepackage{txfonts}
\usepackage{graphicx,color}
\usepackage{epsfig}

\def\ulamek#1#2{\mbox{\normalfont$\frac{#1}{#2}$}}
\def\mmc{\mathcal M}

\title[Multidimensional Catalan and related numbers as Hausdorff moments]{Multidimensional Catalan and related numbers as Hausdorff moments}


\author{Katarzyna G\'{o}rska}                     
\address{H. Niewodnicza\'{n}ski Institute of Nuclear Physics, \\ Polish Academy of Sciences,\\ ul. Eljasza-Radzikowskiego 152, \\ 31342 Krak\'{o}w, Poland}    
\email{katarzyna.gorska@ifj.edu.pl}  

\author{Karol A. Penson}             
\address{Laboratoire de Physique Th\'{e}orique \\ de la Mati\`{e}re Condens\'{e}e (LPTMC), \\ Universit\'{e} Pierre et Marie Curie, \\ CNRS UMR 7600, Tour 13 - 5i\`{e}me \'et., \\Bo\^{i}te Courrier 121, 4 place Jussieu, \\ F 75252 Paris Cedex 05, France}
\email{penson@lptl.jussieu.fr}

\thanks{The authors acknowledge support from Agence Nationale de la Recherche (Paris, France) under Program PHYSCOMB No. ANR-08-BLAN-0243-2 and under the PHC Polonium, Campus France, project no. 288372A}
\begin{document}

\maketitle

\begin{abstract}
We study integral representation of so-called $d$-dimensional Catalan numbers $C_{d}(n)$,  defined by $\left[\prod\limits_{p=0}^{d-1} \ulamek{p!}{(n+p)!}\right] (d n)!$, $d = 2, 3, \ldots$, $n=0, 1, \ldots\,\,$. We prove that the $C_{d}(n)$'s are the $n$th Hausdorff power moments of positive functions $W_{d}(x)$ defined on $x\in[0, d^{\,d}]$. We construct exact and explicit forms of $W_{d}(x)$ and demonstrate that they can be expressed as combinations of $d-1$ hypergeometric functions of type $_{d-1}F_{d-2}$ of argument $x/d^{\,d}$. These solutions are unique. We analyse them analytically and graphically. A combinatorially relevant, specific extension of $C_{d}(n)$ for $d$ even in the form $D_{d}(n)=\left[\prod\limits_{p = 0}^{d-1} \ulamek{p!}{(n+p)!}\right]\left[\prod\limits_{q = 0}^{d/2 - 1} \ulamek{(2 n + 2 q)!}{(2 q)!}\right]$ is analyzed along the same lines.
\end{abstract}


\section{Introduction}

Amongst many existing generalizations of classical Catalan numbers $C(n)=\ulamek{1}{n+1}\binom{2n}{n}$, those that include the parameter that in certain sense can be associated with the spacial dimension $d$, are particularly interesting. They permit to extend to higher dimensions $d > 2$ the notions of objects enumerated by $C(n)$ in $d = 2$. We shall be concerned in this note with one of such generalizations, called $d$-dimensional Catalan numbers \cite{JPAllouche99, RASulanke04, VKlebanov09}, which are defined as:
\begin{equation}\label{e1-1}
C_{d}(n) = \left[\prod_{p = 0}^{d - 1} \frac{p!}{(n+p)!}\right] (d n)!, \quad n = 0, 1, \ldots; \,\, d = 2, 3, \ldots ,
\end{equation}
which for $d = 2$ clearly reduce the conventional Catalan numbers $C(n)$. The form of Eq.~(\ref{e1-1}) guarantees that $C_{d}(0) = 1$ for all $d$. We shall refer to Sloane's Online Encyclopedia of Integer Sequences (OEIS) \cite{wwwOEIS} and quote initial terms, $n=0, 1, \ldots, 7$, of several sequences $C_{d}(n)$, $d = 2, 3, 4$ and $5$, along with the labelling of their entries in the OEIS:
\begin{itemize}
\item for $d = 2$: $1, 1, 2, 5, 14, 42, 132, \ldots$ (A00108), which are the Catalan numbers,

\item for $d = 3$: $1, 1, 5, 42, 462, 6006, 87516, \ldots$ (A005789, A151334),

\item for $d = 4$: $1, 1, 14, 462, 24024, 1662804, 140229804, \ldots$ (A005790), 

\item for $d = 5$: $1, 1, 42, 6006, 1662804, 701149020, 396499770810, \ldots$ (A005791),

\item for general $d$, see A060854.
\end{itemize}

The explicit form of $C_{d}(n)$'s permits one to immediately write down some of their characteristics. If $\Delta(k, a) = \ulamek{a}{k}, \ulamek{a+1}{k}, \ldots, \ulamek{a+k-1}{k}$ denotes a special list of $k$ elements, then the ordinary generating function (ogf) of $C_{d}(n)$'s can be written as:
\begin{equation}\label{e1-2}
g(d, z) = \sum_{n=0}^{\infty} C_{d}(n)\,z^n = {_{d}F_{d-1}}\left({\Delta(d, 1) \atop 2, 3, \ldots, d} \Big\vert d^{\,d} z\right).
\end{equation}
Similarly, the exponential generating function (egf) of $C_{d}(n)$'s reads:
\begin{equation}\label{e1-3}
G(d, z) = \sum_{n=0}^{\infty} C_{d}(n)\,\frac{z^n}{n!} = {_{d}F_{d}}\left({\Delta(d, 1) \atop 1, 2, \ldots, d} \Big\vert d^{\,d} z\right).
\end{equation}
The use of Stirling's formula gives the leading term of $n\to\infty$ asymptotics for $C_{d}(n)$:
\begin{equation}\label{e1-4}
C_{d}(n) \,\,_{\overrightarrow{n\to \infty}}\,\, n^{-\ulamek{d^2 - 1}{2}} d^{\, d n} + \ldots, \quad d = 2, 3, \ldots\,.
\end{equation}
In Eqs.~(\ref{e1-2}) and (\ref{e1-3}) we have employed the standard notation for the generalized hypergeometric function ${_{p}F_{q}}\left({(\alpha_{p}) \atop (\beta_{q})} \Big\vert x \right)$, with $(\alpha_{p})$ and $(\beta_{q})$ the lists of $p$ "upper" and $q$ "lower" parameters, respectively, see \cite{APPrudniko98v3}. Observe that since in ${_{d}F_{d}}$ of Eq.~(\ref{e1-3}) there is a pair of lower and upper parameters differing by one, the appropriate function ${_{d}F_{d}}$ can be reduced to a combination of ${_{d-1}F_{d-1}}$'s, see the formula 7.2.13.17 on page 439 of \cite{APPrudniko98v3}.

Inspired by the very fruitful interpretation of Catalan numbers $C(n)$ as moments of a positive function on $x\in[0, 4]$, which is intimately related to the famous Wigner's semicircle law \cite{LArnold71}, we set out to consider the sequences $C_{d}(n)$, $d > 2$ as Hausdorff power moments and have defined an objective of obtaining for $d > 2$ the equivalents of the solution for $d = 2$, quoted in Eq.~(\ref{e2-8}) below.

The paper is organised as follows: in Sec. 2 we describe the method of obtaining exact and explicit solutions for $d \geq 2$. Subsequently we write down the general solution for $d$ arbitrary and quote the specific cases of $d = 2, 3, 4$ and $5$. In Section 3 we discuss some possible generalizations of $C_{d}(n)$'s. In Sec. 4 we close the note with short conclusions and comments about possible applications of the probability distributions found here.


\section{Solutions of the Hausdorff moment problem}


We are seeking the solutions of the following Hausdorff moment problem:
\begin{equation}\label{e2-1}
\int_{0}^{R(d)} x^n W_{d}(x) \, dx = C_{d}(n), \quad n=0, 1, \ldots, \quad d = 2, 3, \ldots,
\end{equation}
where $R(d)$ - the upper edge of the support of $W_{d}(x)$ - will be determined below. The conventional estimate $R(d)=\lim\limits_{n\to\infty} [C_{d}(n)]^{1/n} = d^{d}$ will be confirmed later by the Mellin transform analysis. As a preliminary step we shall demonstrate that the sought for $W_{d}(x)$ defined in Eq.~(\ref{e2-1}) is positive. By using the Gauss-Legendre multiplication formula for gamma function to Eq.~(\ref{e1-1}) and introducing complex $s$ such that $n = s-1$ we obtain
\begin{equation}\label{e2-2}
C_{d}(s-1) = (2\pi)^{\frac{1-d}{2}} d^{\frac{1}{2}-d} \left(\prod_{k=0}^{d-1} k!\,\right)(d^{\,d})^{\,s}\, \prod_{j=0}^{d-1} \frac{\Gamma\Big(s -1 + \ulamek{j+1}{d}\Big)}{\Gamma(s+j)},
\end{equation} 
which should be interpreted as the Mellin transform of $W_{d}(x)$ i. e. $\int_{0}^{\infty}x^{s-1} W_{d}(x) dx$, denoted by $\mmc\left[W_{d}(x); s\right]$, see Ref.~\cite{IASneddon72}. Since $j > \ulamek{j+1-d}{d}$ for all $0 \leq j \leq d-1$, the individual term labelled by $j$ in the second product of Eq.~(\ref{e2-2}) has the inverse Mellin transform \cite{IASneddon72}, see the formula 8.4.2.3 on page 631 from \cite{APPrudniko98v3}
\begin{equation}\label{e2-3}
\mmc^{-1}\left[\frac{\Gamma\Big(s -1 + \ulamek{j+1}{d}\Big)}{\Gamma(s+j)}; x\right] = \frac{x^{(j+1)/d - 1}\, (1-x)^{j-(j+1)/d}}{\Gamma\Big(1 + j - \ulamek{j+1}{d}\Big)},
\end{equation}
$j=0,1,\ldots d-1$, e. g. it is proportional to the standard probabilistic beta distribution \cite{ChKleiber03} in the variable $x$, which is a positive and absolutely continuous  function for $0 \leq x \leq 1$. We perceive now $C_{d}(s-1)$ as a product of $d$ such individual terms. Then the weight $W_{d}(x)$ is a positive and absolutely continuous function on $[0, R(d)]$, since it is a $d$-fold Mellin convolution of positive and absolutely continuous  functions on $[0, 1]$. In the final result we accommodate the prefactor $(d^{\,d})^s$ which indicates, via elementary property of the Mellin transform \cite{IASneddon72}, that the solution of Eq.~(\ref{e2-1}) will depend on $x/d^{\,d}$.

It turns out that such a $d$-fold Mellin convolution can be carried out explicitly. The key step is first to identify the weight $W_{d}(x)$ as a special case of Meijer $G$ function $G^{m, n}_{p, q}$ \cite{APPrudniko98v3}. This is a direct consequence of \eqref{e2-2} and reads
\begin{equation}\label{e2-4}
W_{d}(x) = (2\pi)^{\frac{1-d}{2}} d^{\frac{1}{2}-d} \left(\prod_{k=0}^{d-1} k!\right)\, G^{d, 0}_{d, d}\left(\frac{x}{d^{d}} \Big\vert {0, 1, \ldots, d-1 \atop-\Delta(d, 0)}\right),
\end{equation}
where $\Delta(n, a) = \ulamek{a}{n}, \ulamek{a+1}{n}, \ldots, \ulamek{a+n-1}{n}$. Next, the Meijer $G$ function is converted to the hypergeometric form, using formulas 16.17.2 and 17.17.3 of \cite{NIST}, which is the Slater theorem. We quote only the final result which is equal to: 
\begin{equation}\label{e2-5}
W_{d}(x) = \sum_{j=1}^{d-1}\, \frac{c_{j}(d)}{x^{j/d}} \, _{d-1}F_{d-2}\left({-\ulamek{j}{d}, -1-\ulamek{j}{d}, \ldots, -d+2-\ulamek{j}{d} \atop
1-\ulamek{1}{d}, 1-\ulamek{2}{d}, \ldots, 1 - \ulamek{j-1}{d}\,;\,\, 1+\ulamek{1}{d}, 1+\ulamek{2}{d}, \ldots, 1+\ulamek{d-j-1}{d}} \Big\vert \frac{x}{d^{\,d}} \right),
\end{equation}
defined for $0 \leq x \leq d^{\,d}$, which implies $R(d) = d^{d}$ in Eq. \eqref{e2-1}. (For the reader's convenience we point out that in Eq.~(\ref{e2-5}), in the lower list of parameters of ${}_{d-1}F_{d-2}$, there are \textit{two} sequences of numbers, which contain $j-1$ and $d-1-j$ terms, respectively). The numerical coefficient $c_{j}(d)$ is equal to
\begin{equation}\label{e2-6}
c_{j}(d) = (2\pi)^{\ulamek{1-d}{2}} d^{j-d+\ulamek{1}{2}} \left[\prod\limits_{p=1}^{d-1} \frac{p!}{\Gamma\Big(p+\ulamek{j}{d}\Big)}\right] \left[\prod\limits_{k=1}^{j-1} \Gamma\Big(\ulamek{k}{d}\Big)\right] \left[\prod\limits_{k=j+1}^{d-1} \Gamma\Big(\ulamek{j-k}{d}\Big)\right],
\end{equation} 
where  $j = 1, \ldots, d-1$ and $d = 2, 3, \ldots\,\,$.

The structure of parameter list of Meijer G function in Eq. \eqref{e2-4} warrants that the assumptions of formula 2.24.2.1 in \cite{APPrudniko98v3} are satisfied:
\begin{equation}\label{e2-7}
-\frac{1}{d} \sum_{k=0}^{d-1} k - \sum_{k=0}^{d-1} k = - \frac{d^{2}-1}{2} < 0, \quad d = 2, 3, \ldots\,.
\end{equation}
Therefore the Mellin transform of $W_{d}(x)$ is well defined for $\Re(s) > \frac{d-1}{d}$.

\begin{figure}[h!]
\includegraphics[scale=0.4]{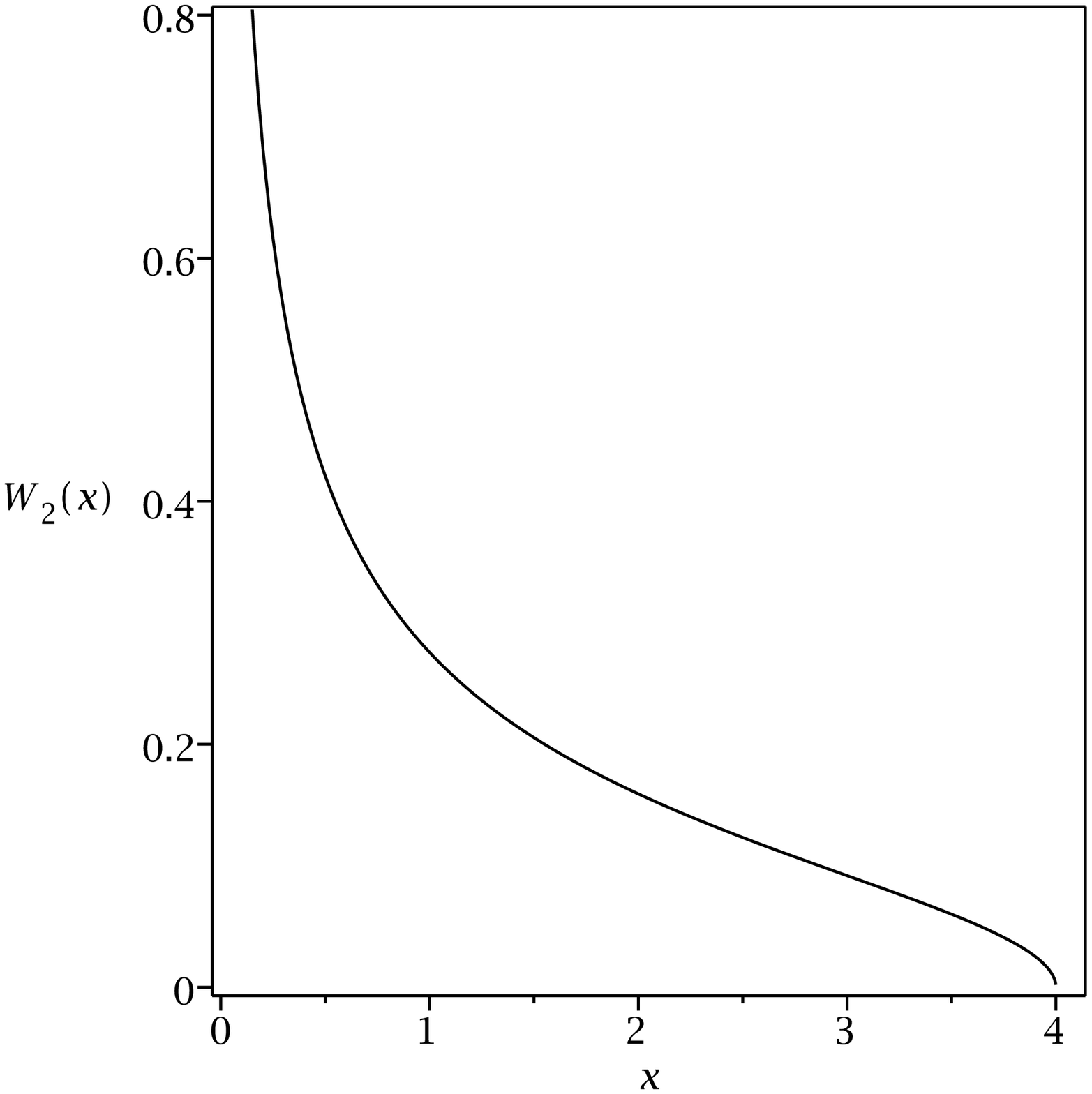}
\caption{\label{fig1} The density $W_{2}(x)$, s. Eq. \eqref{e2-8}.}
\end{figure}

\begin{figure}[h!]
\includegraphics[scale=0.4]{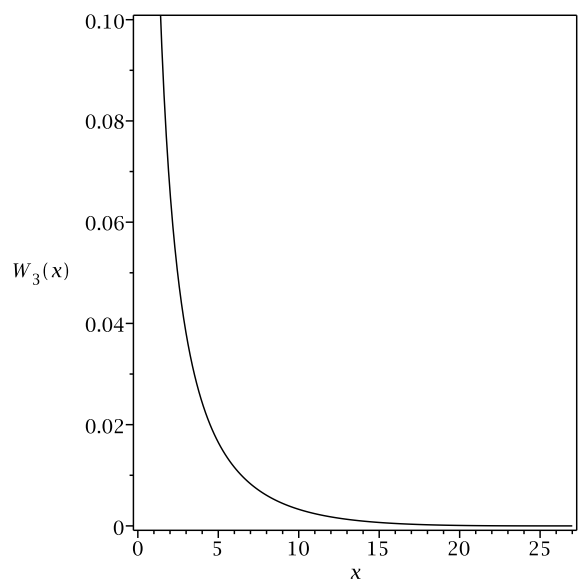}
\caption{\label{fig2} The density $W_{3}(x)$, s. Eq. \eqref{e2-9}.}
\end{figure}

We shall explicitly write down the solutions for $d = 2, 3, 4$ and $5$, starting with $W_{2}(x)$:
\begin{equation}\label{e2-8}
W_{2}(x) = \frac{1}{2\pi}\, \sqrt{\frac{4 - x}{x}}, \quad 0 < x \leqslant 4,
\end{equation}
which is obtained in many references \cite{VAMarchenko67, KAPenson01}, see Fig.~\ref{fig1}. It is the only density that can be expressed by an elementary function. Furthermore, for $d > 2$ no density can be expressed by \textit{standard} special functions, and the hypergeometric form is the final one. For $d~=~3,~\ldots,~5$ the solutions read:
\begin{equation}\label{e2-9}
W_{3}(x) = \frac{c_{1}(3)}{x^{1/3}} \,_{2}F_{1} \left({-\ulamek{4}{3}, -\ulamek{1}{3} \atop \ulamek{4}{3}} \Big\vert \frac{x}{3^3}\right) + \frac{c_{2}(3)}{x^{2/3}} \,_{2}F_{1} \left({-\ulamek{5}{3}, -\ulamek{2}{3} \atop \ulamek{2}{3}} \Big\vert \frac{x}{3^3} \right), \quad 0 < x \leqslant 3^{3},
\end{equation}
\begin{align}\label{e2-10}
W_{4}(x) &= \frac{c_{1}(4)}{x^{1/4}} \,_{3}F_{2}\left({-\ulamek{9}{4}, -\ulamek{5}{4}, -\ulamek{1}{4} \atop \ulamek{5}{4}, \ulamek{3}{2}} \Big\vert \frac{x}{4^4}\right) + \frac{c_{2}(4)}{x^{1/2}} \,_{3}F_{2} \left({-\ulamek{5}{2}, -\ulamek{3}{2}, -\ulamek{1}{2} \atop \ulamek{3}{4}, \ulamek{5}{4}} \Big\vert \frac{x}{4^4}\right) \\[0.7\baselineskip] 
&+ \frac{c_{3}(4)}{x^{3/4}} \,_{3}F_{2}\left({-\ulamek{11}{4}, -\ulamek{7}{4}, -\ulamek{3}{4} \atop \ulamek{1}{2}, \ulamek{3}{4}} \Big\vert \frac{x}{4^4}\right), \quad 0 < x \leqslant 4^{4},\nonumber  
\end{align}
\begin{align}\label{e2-11}
W_{5}(x) &= \frac{c_{1}(5)}{x^{1/5}} \,_{4}F_{3} \left({-\ulamek{16}{5}, -\ulamek{11}{5}, -\ulamek{6}{5}, -\ulamek{1}{5} \atop \ulamek{6}{5}, \ulamek{7}{5}, \ulamek{8}{5}} \Big\vert \frac{x}{5^5}\right) + \frac{c_{2}(5)}{x^{2/5}} \,_{4}F_{3} \left({-\ulamek{17}{5}, -\ulamek{12}{5}, -\ulamek{7}{5}, -\ulamek{2}{5} \atop \ulamek{4}{5}, \ulamek{6}{5}, \ulamek{7}{5}} \Big\vert \frac{x}{5^5}\right)  \\[0.7\baselineskip] 
&+ \frac{c_{3}(5)}{x^{3/5}} \,_{4}F_{3} \left({-\ulamek{18}{5}, -\ulamek{13}{5}, -\ulamek{8}{5}, -\ulamek{3}{5} \atop \ulamek{3}{5}, \ulamek{4}{5}, \ulamek{6}{5}} \Big\vert \frac{x}{5^5}\right) + \frac{c_{4}(5)}{x^{4/5}} \,_{4}F_{3} \left({-\ulamek{19}{5}, -\ulamek{14}{5}, -\ulamek{9}{5}, -\ulamek{4}{5} \atop \ulamek{2}{5}, \ulamek{3}{5}, \ulamek{4}{5}} \Big\vert \frac{x}{5^5}\right), \nonumber \\ 
& 0 < x \leqslant 5^{5}. \nonumber
\end{align}
The coefficients $c_{j}(d)$, $j = 1, \ldots, d-1$, for $d = 3, \ldots, 6$, are collected in Tab.~\ref{tab1}. With $c_{j}(6)$'s given there and using Eqs.~(\ref{e2-5}) and (\ref{e2-6}), the reader can easily reconstruct $W_{6}(x)$, which will not be reproduced here. The solution $W_{3}(x)$ is represented in Fig.~\ref{fig2}.

\begin{table}
\begin{tabular}{c | c c c c c}\hline
j & 1 & 2 & 3 & 4 \\ \hline
$c_{j}(3)$ & $-\frac{3^{3}\sqrt{3}}{16\pi^3} \Gamma\Big(\ulamek{2}{3}\Big)^3$ & $\frac{3^{2}}{10}\Gamma\Big(\ulamek{2}{3}\Big)^{-3}$ & --- & --- & --- \\
$c_{j}(4)$ & $\frac{4^{4} 2}{75 \pi^4} \Gamma\Big(\ulamek{3}{4}\Big)^4$ & $-\frac{4^{3}}{15 \pi^2}$ & $\frac{4^{6}}{4851}\Gamma\Big(\ulamek{3}{4}\Big)^{-4}$ & --- & ---\\ 
$c_{j}(5)$ & $-\frac{5^{9} \sqrt{5} \Gamma\Big(\ulamek{4}{5}\Big)^5\, A^4}{2^{5} 3^{2} 11^{2} \pi^5\, B}$ & $\frac{5^{9}\sqrt{5} \Gamma\Big(\ulamek{3}{5}\Big)^5 B^4}{2^{6} 7^{3} 17 \pi^5\, A}$ & $-\frac{5^{9} \sqrt{5} (AB)^{-1}}{2^{7} 3^{4} 13^{2} \Gamma\Big(\ulamek{3}{5}\Big)^5}$ & $\frac{5^{9} \sqrt{5} (AB)^{-1}}{2^{7} 3^{4} 7^{2} 19 \Gamma\Big(\ulamek{4}{5}\Big)^4}$ & ---\\
$c_{j}(6)$ & $\frac{2^{13} 3^{16}\, E}{7^{4}13^{3} 1805 \pi^6}$ & $-\frac{3^{19} D}{2^{10} 7^{3} 65 \pi^6}$ & $\frac{2^{21}}{7^{2}5^{2}3^{4} \pi^3}$ & $-\frac{3^{17}}{11^{2} 5^{3} 2^{8} 7 D} \,\,$ & $\frac{2^{22} 3^{17}}{11^{4} 5^{4} 17^{3} 23^{2} 29 E}$ 
\end{tabular}
\vspace{6mm}
\caption{\label{tab1} The coefficients $c_{j}(d)$, see Eqs.~(\ref{e2-5}) and (\ref{e2-6}) for $d = 3, \ldots, 6$ and $j = 1, \ldots, d-1$; To simplify the notation in $c_{j}(5)$ and $c_{j}(6)$ we set $A=\sin\Big(\ulamek{\pi}{5}\Big)$, $B=\sin\Big(\ulamek{2\pi}{5}\Big)$ and $E = \Gamma\Big(\ulamek{5}{6}\Big)^6$, $D = \Gamma\Big(\ulamek{2}{3}\Big)^6$.}
\end{table}


\section{Generalization of multidimensional Catalan numbers}

In this paragraph we analyse the extension of $C_{4}(n)$ obtained by replacing $(4n)!$ in Eq.~(\ref{e1-1}) by $(2n)! (2n + 2)!\,$. The corresponding sequence 
$$D_{4}(n)\equiv 6(2n)!(2n + 2)! \left[\prod\limits_{r=0}^{3}(n+r)!\right]^{-1}$$
has attracted attention in several contexts, as it appears in \cite{JPAllouche99, MDSainteCatherine86, OBernardi07, MBousquetMelov08}.

\begin{figure}[!h]
\includegraphics[scale=0.4]{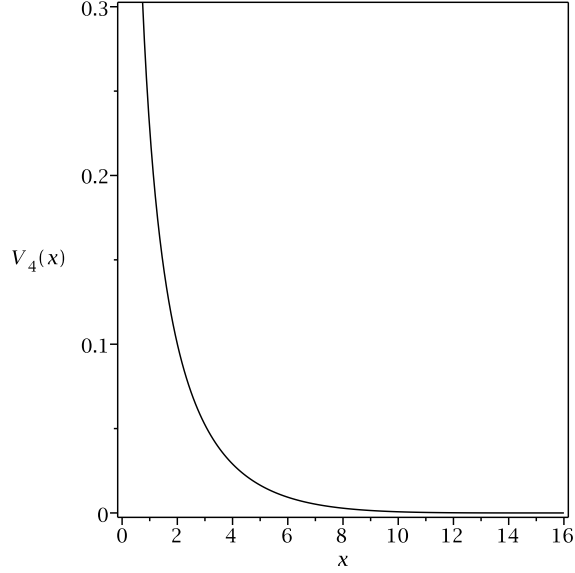}
\caption{\label{fig3} The density $V_{4}(x)$, s. Eq. \eqref{e3-3}.}
\end{figure}

The initial terms of $D_{4}(n)$ are $1, 1, 4, 30, 330, 4719, 81796, 1643356$, for \\$n~=~0, 1, \ldots, 7$. It is listed as A006149 in OEIS where also additional information can be found. It turns out that the ogf of $D_{4}(n)$ can be expressed by the elliptic functions $\mathbb{E}(y)$ and~$\mathbb{K}(y)$ \cite{NIST}:
\begin{align}\label{e3-1}
\sum_{n=0}^{\infty} D_{4}(n)\, z^n &= \frac{1 + 6z}{4 z^2} + \frac{(1 - 16z)\,(1 + 112z)}{240\,\pi\, z^3}\, \mathbb{K}(4\sqrt{z})  \\ \nonumber
&- \frac{1 + 224z + 256 z^2}{240\,\pi\, z^3}\, \mathbb{E}(4\sqrt{z}).
\end{align}
In fact the sequence $D_{4}(n)$ allows for the same kind of analysis as does the ensemble of $C_{d}(n)$'s. The Hausdorff moment problem for $D_{4}(n)$, namely
\begin{equation}\label{e3-2}
\int_{0}^{h} x^n\, V_{4}(x) \, dx = D_{4}(n) = \frac{6 (2n)! (2n + 2)!}{\prod\limits_{r=0}^{3} (n+r)!}, \quad n = 0, 1, \ldots,
\end{equation}
can be solved by the method of Mellin convolution and the use of Meijer G function elucidated above. The weight can be proven to be positive on $x\in[0, h]$ with $h~=~16$ and reads:
\begin{align}\label{e3-3}
V_{4}(x) &= \ulamek{1}{15 \pi^2} \left[\left(\ulamek{64}{\sqrt{x}} + 56\sqrt{x} + \ulamek{x^{3/2}}{4}\right) \, \mathbb{E}\left(\sqrt{1 - \ulamek{x}{16}}\right) - 2\sqrt{x} (16 + x) \, \mathbb{K}\left(\sqrt{1 - \ulamek{x}{16}}\right)\right]. 
\end{align}
We plot the function $V_{4}(x)$ on Fig. \ref{fig3}.

\begin{figure}[!h]
\includegraphics[scale=0.4]{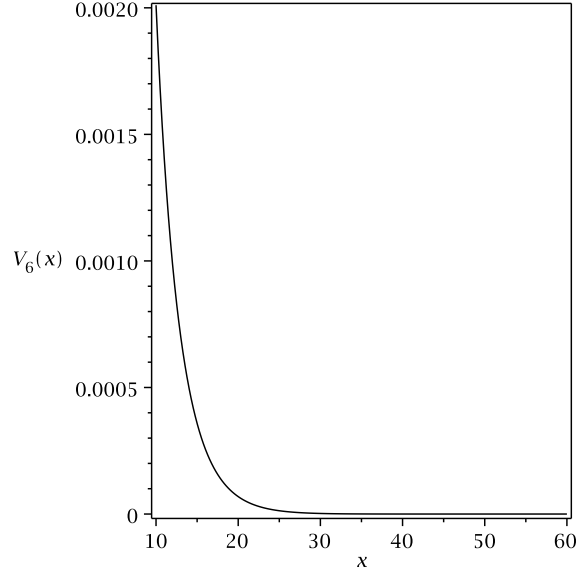}
\caption{\label{fig4} The density $V_{6}(x)$, s. Eq. \eqref{e3-8} for $d=6$.}
\end{figure}

The sequence $D_{4}(n)$ analysed above is a special case $d=4$ of the following generalization of $C_{d}(n)$ defined for \textit{even} $d$:
\begin{equation}\label{e3-4}
D_{d}(n) = \left[\prod_{r = 0}^{d-1} \frac{r!}{(n+r)!}\right] \, \left[\prod_{s = 0}^{d/2 - 1} \frac{(2 n + 2 s)!}{(2 s)!}\right], \qquad d = 2, 4, 6, \ldots,
\end{equation}
in which the parameter $d$ should not be associated anymore with the spacial dimension. Several exact characteristics of the sequences $D_{d}(n)$ are available. The ordinary generating function reads
\begin{equation}\label{e3-5}
\tilde{g}(d, z) = {_{d/2}F_{d/2-1}}\left({\ulamek{1}{2}, \ulamek{3}{2}, \ldots, \ulamek{d-1}{2}, \, 1 \atop \ulamek{d}{2} + 1, \ulamek{d}{2} + 2, \ldots, d} \Big\vert 2^{d} z \right),
\end{equation}
whereas the corresponding exponential generating function is equal to
\begin{equation}\label{e3-6}
\tilde{G}(d, z) = {_{d}F_{d}}\left({\ulamek{1}{2}, \ulamek{3}{2}, \ldots, \ulamek{d-1}{2}, \, 1, 2, \ldots, \ulamek{d}{2} \atop 1, 2, 3, \ldots, d} \Big\vert 2^{d} z \right).
\end{equation}
The leading term of the $n\to\infty$ asymptotics for $D_{d}(n)$ can be obtained by using the Stirling's formula and it has the following form
\begin{equation}\label{e3-7}
D_{d}(n)  \,\,_{\overrightarrow{n\to \infty}}\,\, n^{-d\,(d - 1)/4} \, 2^{d\, n}, \quad d=4, 6, \ldots\, .
\end{equation}

It is remarkable that the Hausdorff moment problem for $D_{d}(n)$, i. e.
\begin{equation*}
\int_{0}^{\kappa(d)} x^n\, V_{d}(x)\, dx \,=\, D_{d}(n), \quad n=0, 1, \ldots; \quad d= 4, 6, \ldots
\end{equation*}
can be exactly solved as well in terms of positive functions $V_{d}(x)$ defined on $x~\in~[0, 2^{d}]$, i. e. $\kappa(d) = 2^{d}$, which read:
\begin{align}\label{e3-8}
V_{d}(x) &= \frac{2^{-d}\, \prod\limits_{r=0}^{d-1} r!}{\prod\limits_{k = 0}^{d/2 - 1} \Gamma\Big(k + \ulamek{1}{2}\Big)\, k!} \, G^{\,d/2,\, 0}_{d/2,\, d/2}\left(\frac{x}{2^d} \Big\vert {\ulamek{d}{2}, \ulamek{d}{2} + 1, \ldots, d-1 \atop -\ulamek{1}{2}, \ulamek{1}{2}, \ulamek{3}{2}, \ldots, \ulamek{d-3}{2}}\right).
\end{align}
Here the condition 2.24.2.1 in \cite{APPrudniko98v3} implies $-\ulamek{d}{4} (d+1) < 0$, $d=4, 6, \ldots$, which is always satisfied. In addition, the Mellin transform of $V_{d}(x)$ is well defined for $\Re(s) > \ulamek{1}{2}$. The proof of positivity of $V_{d}(x)$ can be carried out along the lines exposed in the previous Section.

Since in Eqs.~(\ref{e3-8}) in both parameter lists in the Meijer G function there are index pairs that differ by an integer, this Meijer G function cannot be represented by a sum of generalized hypergeometric functions. However the expression \eqref{e3-8} can be easily manipulated algebraically and represented graphically \cite{Map}. In Fig. \ref{fig4} we display $V_{6}(x)$ in the range $x\in[10,\, 60]$. Observe the rapid decrease of this function for $x \gtrsim 25$, followed by a large region where it is practically flat and equals to zero. Similar behavior is observed for higher values of $d$.

\section{Discussion and Conclusions}

We have treated in this work essentially two generalizations of conventional Catalan numbers, which are related to such notions as Young tableaux, hook lengths, generalized Dyck paths, etc. \cite{EDeutsch}. They all turn out to be moments of positive functions on supports included in the positive half line. The relevant weight functions have been obtained explicitly and analyzed graphically. All these positive functions are unique solutions of Hausdorff moment problems. The key tools in this approach had been the inverse Mellin transform and the encoding with Meijer $G$ functions. The positivity of solutions has been rigorously proven using the method of Mellin convolution, applied to related problems previously \cite{KGorska10, KAPenson11, WMlotkowski12}.

It needs to be specified that the function $W_{2}(x)$ of Eq. \eqref{e2-8} is the known Marchenko-Pastur distribution \cite{VAMarchenko67, KAPenson11} which describes the level statistic of random Wishart matrices $W = G G^{\dag}$, where $G$ is a square, $N\times N$ random Ginibre matrix. As far as applications for random matrices are concerned two problems appear to be relevant for the distributions found in the present work. 

First it would be intriguing to know if the distributions $W_{d}(x)$ for $d \geq 3$, and $V_{d}(x)$ for $d = 4, 6, \ldots$ would correspond to limit spectral densities of certain (if any) ensembles of random matrices. A second possibility is to extend the analysis of products of square random matrices to products of rectangular $N\times M$ random matrices with $r = N/M$. A case in point is a detailed analysis of products of rectangular Gaussian random matrices carried out in \cite{ZBurda10}. Therefore, once the relevant matrix ensemble has been properly identified, it is quite feasible to undertake the analysis of appropriate products of rectangular matrices. This would lead, in the spirit of \cite{ZBurda10} to, for instance, $W^{(r)}_{3}(x)$ parametrized by $r$, with $W^{(1)}_{3}(x) \equiv W_{3}(x)$. Both of these problems are under active consideration.

\section{Acknowledgement}

We thank G.~H.~E.~Duchamp and Olivier G\'{e}rard for useful suggestions and discussions. 
\bibliographystyle{amsplain}

\begin{thebibliography}{99}
\bibitem{Map} We have made extensive use of Maple$^{\tiny{\textregistered}}$ in this work.
\end{thebibliography}

\end{document}